\documentclass{elsarticle}
\usepackage{amssymb,wrapfig,amsmath, amsthm, amscd,ifthen}
\usepackage[cp1251]{inputenc}
\usepackage{color} 
\usepackage{hyperref} 
\usepackage{graphicx,caption,subcaption,mathrsfs,appendix}
\newtheorem{theorem}{Theorem}

\newcommand{\F}{\mathcal{F}}

\newcommand{\E}{\mathbb E}
\newcommand{\N}{\mathbb N}
\newcommand{\1}{\mathbf 1}

\begin{document}
\title{Preferential attachment with fitness dependent choice}
\tnotetext[t1]{The present work was funded by a grant from the Russian Science Foundation (project No. 19-71-00043).}
\author[1]{Y.A. Malyshkin}
\ead{yury.malyshkin@mail.ru}

\address[1]{Tver State University//Moscow Institute of Physics and Technology}

\begin{abstract}
We study the asymptotic behavior of the maximum degree in the evolving tree model with a choice based on both degree and fitness of a vertex. The tree is constructed in the following recursive way. Each vertex is assigned a parameter to it that is called a fitness of a vertex. We start from two vertices and an edge between them. On each step we consider a sample with repetition of $d$ vertices, chosen with probabilities proportional to their degrees plus some parameter $\beta>-1$. Then we add a new vertex and draw an edge from it to the vertex from the sample with the highest product of fitness and degree. We prove that dependent on parameters of the model, the maximum degree could exhibit three types of asymptotic behavior: sublinear, linear and of $n/\ln n$ order, where $n$ is the number of edges in the graph.  
\end{abstract}


\maketitle

\section{Introduction}

In the present work, we study the combination of two modification of the standard preferential attachment graph model, which was introduced in \cite{BA99}. The preferential attachment graph is constructed in the following way. First, we start with some initial graph $G_1$, usually, for simplification purposes, it consists of two vertices and an edge between them. Then on each step we add a new vertex and draw an edge from it to an already existing vertex, that is chosen by some rule. For the standard preferential attachment model, the rule is that we choose a vertex with probability proportional to its degree. Usually, the one considers the rule when we choose a vertex with probability proportional to its degree plus some parameter $\beta>-1$ (see, e.g. \cite{Mori02,Mori05}). Such a model was widely studied (see, e.g., \cite{Hof16}, section 8) and different modification have been introduced.

In the present work, we study the combination of two modifications of this model. One of the modification is the introduction of the fitness to the model (see, e.g., \cite{BCDR07}). Fitness is a parameter that assigns to each vertex and affects its probability to be chosen at each step. The other modification is the introduction of a choice to the model (see, e.g., \cite{KR14,MP14,MP15}). In this modification, we consider the sample of $d$ independently chosen vertices and then choose one of them by some rule. Two types of rules have been considered: degree base rule (see, e.g., \cite{HJ16,Mal18}) and location (or fitness) based choice (see, e.g., \cite{GLY19,HJY20}). In the present work, we consider choice based on both degree and fitness.

Let us introduce our model. Fix $\beta>-1$ and $d\in\N$, $d>1$. We consider a sequence of graphs $G_n$ build recursively as following. We start with initial graph $G_1$ that consists of a two vertex $v_0$, $v_1$ and an edge between them. Graph $G_{n+1}$ is built from $G_n$ by adding a new vertex $v_{n+1}$ and drawing an edge from it to the vertex of $G_n$, choosen by following rule. We first consider a sample of $d$ vertices of $G_n$ choosen independently from each other with probabilities proportional to their degrees plus parameter $\beta$. Then we choose vertex among them that maximize function $W(v,n):=\lambda_i\deg_{G_n} v$
where $\lambda_i$ is the fitness of $v_i$. We consider case when $\lambda_i$ are i.i.d. random variables that take two values, $1$ and $\lambda>1$ with non-zero probabilities. For simplification we also suggest that $\lambda$ is not rational so $W(v,n)$ would not take the same value on vertices with different fitness.

Let us formulate our main result. Let $M(n)$ be the maximum degree of vertices of $G_n$.

\begin{theorem}
In defined model,
\begin{enumerate}
\item If $d<2+\beta$ than for any $\epsilon>0$ 
$$\Pr\left(\forall n>n_0: n^{\frac{d}{2+\beta}-\epsilon}<M(n)<n^{\frac{d}{2+\beta}+\epsilon}\right)\to 1$$
as $n_0\to\infty$.
\item If $d=2+\beta$, than almost surely
$$\liminf_{n\to\infty}\frac{M(n)\ln n}{n}\geq \frac{2d}{(d-1)\lambda},$$
$$\limsup_{n\to\infty}\frac{M(n)\ln n}{n}\leq \frac{2d}{d-1}.$$
\item If $d>2+\beta$, than almost surely
$$\liminf_{n\to\infty}\frac{M(n)}{n}\geq \frac{x^{\ast}}{\lambda},$$
$$\limsup_{n\to\infty}\frac{M(n)\ln n}{n}\leq x^{\ast},$$
where $x^{\ast}$ is a unique positive root of equation $1-\left(1-\frac{x}{2+\beta}\right)^d=x$.
\end{enumerate}
\end{theorem}
The result is similar to the Theorem 1.1 of \cite{Mal18} and shows that the addition of fitness to the choice from the sample does not affect the type of asymptotic of the maximal degree. Let us provide an outline of the proof. For each case, the lower and the upper bound is proven separately. The general idea is that we obtain either a lower or upper bound for the conditional probability to increase $M(n)$ (or its modification for the lower bound) as a function of $M(n)$. Then we study the properties of this function to get an estimate for $M(n)$. One of the key factors is the coefficient of the first term of the expansion by degrees of $\frac{M(n)}{n}$, which equals to $\frac{d}{2+\beta}$. In case $d<2+\beta$ we would analyse the fraction $\frac{M(n+1)}{M(n)}$ to prove sublinear behavior of $M(n)$. In case $d=2+\beta$ we would construct additional expressions to outline the second term of the expansion to get $\frac{n}{\ln n}$ bounds and for $d>2+\beta$, we would use a stochastic approximation to get linear estimates. 

Let us give a short description of the stochastic approximation approach (see, e.g., \cite{Chen03,Pem07} for more details). Process $Z(n)$ is a stochastic approximation process if it could be written as
$$Z(n+1)-Z(n)=\gamma_n\left(F(Z(n))+E_n+R_n\right)$$
where $\gamma_n$, $E_n$ and $R_n$ satisfy the following condition. $\gamma_n$ is not random and $\sum_{n=1}^{\infty}\gamma_n>\infty$, $\sum_{n=1}^{\infty}(\gamma_n)^2<\infty$, usually one puts $\gamma_n=\frac{1}{n}$ or $\gamma_{n}=\frac{1}{n+1}$. The term $E_n$ is $\F_n$-measurable where $\F_n$ is the natural filtration of $Z(n)$, $\E(E_n|\F_n)=0$ and $\E ((E_n)^2|\F_n)<c$ for some fixed constant $c$. We cosider $E_n=\frac{1}{\gamma_n}(Z(n+1)-\E(Z(n+1)|\F_n))$ and therefore the function $F(x)$ could be found from representation $\E(Z(n+1)-Z(n)|\F_n)=\gamma_n(F(Z(n))+R_n)$ where $R_n$ is a small error term that satisfies $\sum_{n=1}^{\infty}\gamma_n|R_n|<\infty$ almost surely. If these conditions hold (they could be easily checked in our case) then $Z(n)$ converges to the zero set of $F(x)$. In case when there is more than one eligible (nonnegative) root of $F(x)$ we would also prove non-convergence to one of the roots.

The other argument we would use is the persistent hub type of argument (see, e.g., \cite{Gal16} and Proposition 1.2 in \cite{Mal18}). This argument is based
on the following urn model property. If we have random walk $(x(n),y(n))$ that takes steps up and right with probabilities $\frac{x(n)+\beta}{x(n)+y(n)+2\beta}$ and $\frac{y(n)+\beta}{x(n)+y(n)+2\beta}$ (it also represents the evolution of the urn with $x(n)$ white and $y(n)$ black balls in urm model, see e.g., Theorem 3.2 in \cite{M09} or Section 4.2 in \cite{JK77}) then it converges in distribution to continuous Beta-distributions and hence one of the variables $x(n),y(n)$ would exceed the other after some random moment. Also if such process starts at point $(1,a)$, then probability that $x(n)$ exceeds $y(n)$ at some moment would decay exponentially with $a$. Hence, to prove the existence of the persistence hub it is enough to show that the pair $(\deg_{G_n} v_i,\deg_{G_n} v_j)$ dominates urn model in a sense that it has a higher conditional probability to increase the degree of a vertex with a higher degree. We would use such an argument separately for different fitnesses.

\section{Lower bounds}

To prove the lower bounds we separately consider maximum degrees among vertices with different fitnesses and estimate dynamics of the product of fitness and degree.

Let $M_1(n)$ be the maximum degree among vertices with fitness $1$ in $G_n$, $M_{\lambda}(n)$ be the maximum degree among vertices with fitness $\lambda$ in $G_n$.

Let $F_1(k,n)$ be the total weight of vertices with fitness $1$ with degrees more then $k$ in $G_n$,
$F_{\lambda}(k,n)$ be the total weight of vertices with fitness $\lambda$ with degrees more than $k$ in $G_n$. Let $L_1(k,n)$ and $L_{\lambda}(k,n)$ be the number of vertices in $G_n$ with degree $k$ and fitnesses $1$ and $\lambda$ correspondingly. Consider functions 
$$f_n(x,y):=\left(1-\frac{y}{(2+\beta)n}\right)^d-\left(1-\frac{x+y}{(2+\beta)n}\right)^d,$$
$$g_n(x,y):=\sum_{k=0}^{d-1}\left(1-\frac{y}{(2+\beta)n}\right)^{k}\left(1-\frac{x+y}{(2+\beta)n}\right)^{d-1-k}.$$
Note that $f_n(x,y)=\frac{x}{(2+\beta)n}g_n(x,y)$ and both functions are decreasing with $y$ when $x,y\geq 0,$ $x+y\leq (2+\beta)n$. Also, $f_n(x,y)$ is increasing with $x$ while $g_n(x,y)$ is decreasing with $x$, and hence $\lambda f_n(x,y)\geq f_n(\lambda x,y)$.
Then
$$\E\left(M_1(n+1)-M_1(n)|\F_n\right)
=f_n((M_{1}(n)+\beta)L_{1}(M_1(n),n),F_{\lambda}(M_1(n)/\lambda,n)),$$
$$\E\left(M_{\lambda}(n+1)-M_{\lambda}(n)|\F_n\right)
=f_n((M_{\lambda}(n)+\beta)L_{\lambda}(M_{\lambda}(n),n),F_{1}(\lambda M_\lambda(n),n)).$$
Note that either $F_{1}(\lambda M_{\lambda}(n),n)$ or $F_{\lambda}(M_1(n)/\lambda,n)$ equals to $0$. Moreover, if $F_{\lambda}(M_1(n)/\lambda,n)\neq 0$ there is a vertex with fitness $1$ and a degree of at least $\lambda M_{\lambda}(n)$. Let us define process $X_n:=\max\{M_{1}(n),\lambda M_{\lambda}(n)\}$. Note that $M(n)\geq \frac{M_n}{\lambda}$.
Then
$$\E\left(X_{n+1}-X_n|\F_n\right)
\geq\1\{M_{\lambda}(n)>M_{1}(n)/\lambda\}\lambda f_n(M_{\lambda}(n)+\beta,0)$$
$$
+\1\{M_{\lambda}(n)<M_{1}(n)/\lambda\}f_{n}(M_{1}(n)+\beta,0)$$
$$\geq\1\{M_{\lambda}(n)>M_{1}(n)/\lambda\}f_n(\lambda (M_{\lambda}(n)+\beta),0)
$$
$$
+\1\{M_{\lambda}(n)<M_{1}(n)/\lambda\}f_{n}(M_{1}(n)+\beta,0)$$
$$\geq f_{n}(X_n+\beta,0)=\frac{X_n+\beta}{(2+\beta)n}g_{n}(X_n+\beta,0)\geq\frac{X_n}{(2+\beta)n}g_{n}(X_n,0)-\frac{1}{n}.$$
Hence
$$\E\left(\frac{X_{n+1}}{X_{n}}|\F_n\right)\geq 1+\frac{1}{(2+\beta)n}g_{n}(X_n,0)-\frac{1}{nX_n}.$$
Note that $g_{n}(x_n,0)\to d$ if $\frac{x_n}{n}\to 0$. Also, $\prod_{k=1}^{n}\left(1+\frac{d}{(2+\beta)k}\right)$ is of order $n^{\frac{d}{(2+\beta)}}$. Let $\frac{d}{(2+\beta)}\leq 1$. Then, for any $\epsilon$ on the event $\{X_n<n^{\frac{d}{(2+\beta)}-\epsilon}\}$ $X_{n}$ would grow faster than $n^{\frac{d}{(2+\beta)}-\epsilon/2}$ and hence for any $n_0\in\N$ with high probability at some time $n\geq n_0$ $X_n$ would exceed $n^{\frac{d}{(2+\beta)}-\epsilon}$. Also, due to standard large deviation estimates, if $X_{n_0}>(1-\delta)n_0^{\frac{d}{(2+\beta)}-\epsilon}$ then probability that process $X_n$, $n\geq n_0$ would cross a line $(1-2\delta)n^{\frac{d}{(2+\beta)}-\epsilon}$ before it crosses a line $n^{\frac{d}{(2+\beta)}-\epsilon}$ does not exceed $ce^{-n_0^{\frac{d}{(2+\beta)}-\epsilon}}$ for some $c=c(\delta)$. Therefore (for $\frac{d}{(2+\beta)}\leq 1$) with high probability $\liminf_{n\to\infty}\frac{X_n}{n^{\frac{d}{(2+\beta)}-\epsilon}}\geq 1$. 

Similarly, for $\frac{d}{(2+\beta)}>1$ with high probability $\liminf_{n\to\infty}\frac{X_n}{n}>0$. Also, we get that 
$$\E\left(X_{n+1}-X_n|\F_n\right)=f_n(X_n,0)+O\left(\frac{1}{n}\right).$$
Hence, if we define $Z_n=\frac{X_n}{n}$, we would get that
$$\E\left(\left.Z(n+1)-Z(n)\right|\F_n\right)\leq\frac{1}{n+1}\left(\E\left(X_{n+1}-X_n|\F_n\right)-Z_n\right)$$
$$=\frac{1}{n+1}\left(f(Z_n)-Z_n+O\left(\frac{1}{n}\right)\right).$$
Note that if $d>2+\beta$ then due to concavity of $f(x)$ equation $f(x)-x=0$ has two roots in $[0,2+\beta]$ ($0$ and a positive root $x^{\ast}$). Since $Z_n$ does not converge to $0$, by stochastic approximation we get that $\limsup_{n\to\infty} Z_n\leq x^{\ast}$ almost surely, which gives us the lower bound for $d>2+\beta$.

\section{Upper bounds}

To prove the upper bounds we first study dynamic of pairs $(\deg_{G_n}v_i,\deg_{G_n}v_j)$ separately for vertices of fitness $1$ and $\lambda$ to prove the existence of persistence hub among vertices of each fitness. Then we would use it to remove terms $L_{1}(M_{1}(n),n)$ and $L_{\lambda}(M_{\lambda}(n),n)$ and get upper bounds for the increment of the maximum degree and prove upper bounds.

For any vertex $v_i$ with fitness $1$ we get for $n\geq i$
$$\E\left(\deg_{G_{n+1}}(v_i)-\deg_{G_n}(v_i)|\F_n\right)=
$$
$$
=\frac{f_n((\deg_{G_n}(v_i)+\beta)L_{1}(\deg_{G_n}(v_i),n),F_{1}(\deg_{G_n}(v_i),n)+F_{\lambda}(\deg_{G_n}(v_i)/\lambda,n))}{L_{1}(\deg_{G_n}(v_i),n)}
$$
$$
=\frac{\deg_{G_n}(v_i)+\beta}{(2+\beta)n}\times 
$$
$$
\times g_n((\deg_{G_n}(v_i)+\beta)L_{1}(\deg_{G_n}(v_i),n),F_{1}(\deg_{G_n}(v_i),n)+F_{\lambda}(\deg_{G_n}(v_i)/\lambda,n))
$$
and for any vertex $v_i$ with fitness $\lambda$ we get for $n\geq i$
$$\E\left(\deg_{G_{n+1}}(v_i)-\deg_{G_n}(v_i)|\F_n\right)=
$$
$$
=\frac{f_n((\deg_{G_n}(v_i)+\beta)L_{\lambda}(\deg_{G_n}(v_i),n),F_{1}(\lambda\deg_{G_n}(v_i),n)+F_{\lambda}(\deg_{G_n}(v_i),n))}{L_{\lambda}(\deg_{G_n}(v_i),n)}
$$
$$
=\frac{\deg_{G_n}(v_i)+\beta}{(2+\beta)n}\times 
$$
$$
\times g_n((\deg_{G_n}(v_i)+\beta)L_{\lambda}(\deg_{G_n}(v_i),n),F_{1}(\lambda\deg_{G_n}(v_i),n)+F_{\lambda}(\deg_{G_n}(v_i),n)).
$$
For vertices $v_i$ and $v_j$ with the same fitness let us estimate the probability to draw an edge to $v_i$ conditioned on the event that edge is drawn to one of them. Let $\deg_{G_n} v_i>\deg_{G_n} v_j$ (if vertices have the same degree the probability is $1/2$). Note that $g(x,y)$ is decreasing with $y$ and $x+y$. Also, for vertices with fitness $1$ we get that $$F_{1}(\deg_{G_n}(v_i),n)+F_{\lambda}(\deg_{G_n}(v_i)/\lambda,n)\leq F_{1}(\deg_{G_n}(v_j),n)+F_{\lambda}(\deg_{G_n}(v_j)/\lambda,n),$$
$$F_{1}(\deg_{G_n}(v_i),n)+F_{\lambda}(\deg_{G_n}(v_i)/\lambda,n)+(\deg_{G_n}(v_i)+\beta)L_{\lambda}(\deg_{G_n}(v_i),n)
$$
$$\leq F_{1}(\deg_{G_n}(v_j),n) +F_{\lambda}(\deg_{G_n}(v_j)/\lambda,n)+(\deg_{G_n}(v_j)+\beta)L_{\lambda}(\deg_{G_n}(v_j),n).$$
Therefore
$$\Pr\left(\deg_{G_{n+1}}(v_i)-\deg_{G_n}(v_i)=1|\F_n,\right.
$$
$$
\left.\deg_{G_{n+1}}(v_i)-\deg_{G_n}(v_i)+\deg_{G_{n+1}}(v_j)-\deg_{G_n}(v_j)=1\right)$$
$$\geq \frac{\deg_{G_n}(v_i)+\beta}{\deg_{G_n}(v_i)+\deg_{G_n}(v_j)+2\beta}.$$
Hence, after some random time $N$ both $L_{1}(M_{1}(n),n)$ and $L_{\lambda}(M_{\lambda}(n),n)$ would be equal to $1$.

Therefore, for $n>N$ we would get
$$\E\left(M_1(n+1)-M_1(n)|\F_n\right)
=f_n((M_{1}(n)+\beta),F_{\lambda}(M_1(n)/\lambda,n))\leq f_n((M_{1}(n)+\beta),0),$$
$$\E\left(M_{\lambda}(n+1)-M_{\lambda}(n)|\F_n\right)
=f_n((M_{\lambda}(n)+\beta),F_{1}(\lambda M_{\lambda}(n),n))
\leq f_n((M_{\lambda}(n)+\beta),0).$$

Recall that if $d>2+\beta$, then equation $f(x)-x=0$ has two roots in $[0,2+\beta]$ ($0$ and a positive root $x^{\ast}$), and if $d\leq 2+\beta$, then $0$ is the only root in $[0,2+\beta]$.   
Let us define $Z(n):=\frac{M(n)}{n}$. Then
$$\E\left(\left.Z(n+1)-Z(n)\right|\F_n\right)\leq\frac{1}{n+1}\left(f_n(M_n+\beta)-Z_n\right)$$
$$=\frac{1}{n+1}\left(f(Z_n)-Z_n+O\left(\frac{1}{n}\right)\right).$$
For $d\leq 2+\beta$ by stochastic approximation, we get $\limsup_{n\to\infty} Z_n=0$ almost surely. For $d>2+\beta$ by stochastic approximation, we get that $\limsup_{n\to\infty} Z_n\leq x^{\ast}$ almost surely, which gives us an upper bound for $d>2+\beta$.

Now consider the case $d<2+\beta$. Similarly to the lower bound, we get that
$$\E\left(\frac{M(n+1)}{M(n)}|\F_n\right)\leq 1+\frac{f_{n}(M(n)+\beta,0)}{M(n)}=1+\frac{g(M(n)+\beta,0)}{(2+\beta)n}\leq 1+\frac{d}{(2+\beta)n},$$
where in the last inequality we used that $g(x,0)\leq d$ for $x\geq 0$. Hence, as in the lower bound argument, we get that for any $\epsilon>0$ $\limsup_{n\to\infty}\frac{M(n)}{n^{\frac{d}{2+\beta}+\epsilon}}=0$, which gives us the upper bound for $d<2+\beta$.

To get the upper bound for the case $d=2+\beta$, for $c>0$ consider  variables 
$$U_n:=ne^{-c\frac{n}{M(n)}}.$$ 
We get that
$$\E\left(\left.\frac{U_{n+1}}{U_n}\right|\F_n\right) =\frac{n+1}{n}\E\left(\left.e^{c\frac{n}{M(n)}-c\frac{n+1}{M(n+1)}}\right|\F_n\right)$$
$$=\left(1+\frac{1}{n}\right)\E\left(\left.e^{\frac{cn(M(n+1)-M(n))}{X_nX_{n+1}}-\frac{c}{M(n+1)}}\right|\F_n\right)$$
$$=\left(1+\frac{1}{n}\right)\E\left(\left.1-\frac{c}{M(n+1)}+\frac{cn(M(n+1)-M(n))}{M(n)M(n+1)} \right.\right.
$$
$$
\left.\left.+O\left(\left(\frac{1}{M(n+1)}\right)^2+\left(\frac{n}{(M(n+1))^2}\right)^2\right)\right|\F_n\right)$$
$$=\left(1+\frac{1}{n}\right)\E\left(\left.1-\frac{c}{M(n)}+\frac{cn(M(n+1)-M(n))}{(M(n))^2} \right.\right.
$$
$$
\left.\left.+O\left(\left(\frac{1}{M(n)}\right)^2+\left(\frac{n}{(M(n))^2}\right)^2\right)\right|\F_n\right)$$
$$=\left(1+\frac{1}{n}\right)\left(1-\frac{c}{M(n)}+\frac{cn\E\left(\left.M(n+1)-M(n)\right|\F_n\right)}{(M(n))^2} +o\left(\frac{1}{n^{3/2}}\right)\right)$$
$$
\leq\left(1+\frac{1}{n}\right)\left(1-\frac{c}{M(n)}+\frac{cg_{n}(M(n)+\beta,0)}{(2+\beta)M(n)} +o\left(\frac{1}{n^{3/2}}\right)\right)
$$
Note that $g_n(x,0)=d-\frac{d(d-1)}{2}\frac{x}{(2+\beta)n}+o(\frac{x}{n})$ when $\frac{x}{n}\to 0$. Hence if $2+\beta=d$ we get 
$$
\E\left(\left.\frac{U_{n+1}}{U_n}\right|\F_n\right) $$
$$
\leq 1\!+\!\frac{1}{n}\!-\!\frac{c}{M(n)} \!+\!\frac{c}{(2\!+\!\beta)M(n)}\!\left(d\!-\!\frac{d(d\!-\!1)}{2}\frac{M(n)}{(2\!+\!\beta)n}\!+\!o\!\left(\!\left(\frac{M(n)}{2}\right)\!^2\right)\!\right) +o\left(\frac{1}{n}\right)
$$
$$
=\left(1+\frac{1}{n}-\frac{c(d-1)}{2dn}+o\left(\frac{1}{n}\right)\right)
$$
Hence, if $c>\frac{2d}{d-1}$ then $U_n$ is supermartingale, and by Doob's theorem there is a random variable $R$, such that $\sup_n U_n<R$ almost surely. Hence for all (large enought) $n$ 
$$M(n)<\frac{cn}{\ln n - \ln R}$$
almost surely. Therefore 
$$\liminf \frac{\ln n M(n)}{n}\leq \frac{2d}{d-1}$$
almost surely.



\section{Lower bound for $d= 2+\beta$}

Let now prove lower bound for $d= 2+\beta$. The argument is similar to the argument for the upper bound. Recall that $X_n=\max\{M_{1}(n),\lambda M_{\lambda}(n)\}$.

For $\frac{d}{(2+\beta)}=1$ and $c>0$ let consider 
$$Y_n:=e^{c\frac{n}{X_n}}/n$$ 
Note that
$$\E\left(\left.\frac{Y_{n+1}}{Y_n}\right|\F_n\right) =\frac{n}{n+1}\E\left(\left.e^{c\frac{n+1}{X_{n+1}}-c\frac{n}{X_n}}\right|\F_n\right)$$
$$=\frac{n}{n+1}\E\left(\left.e^{\frac{c}{X_{n+1}}-\frac{cn(X_{n+1}-X_n)}{X_nX_{n+1}}}\right|\F_n\right)$$
$$=\frac{n}{n+1}\E\left(\left.1+\frac{c}{X_{n+1}}-\frac{cn(X_{n+1}-X_n)}{X_nX_{n+1}} +O\left(\left(\frac{1}{X_{n+1}}\right)^2+\left(\frac{n}{X_{n+1}^2}\right)^2\right)\right|\F_n\right)$$
$$=\frac{n}{n+1}\E\left(\left.1+\frac{c}{X_{n}}-\frac{cn(X_{n+1}-X_n)}{X_n^2} +O\left(\left(\frac{1}{X_{n}}\right)^2+\left(\frac{n}{X_{n}^2}\right)^2\right)\right|\F_n\right)$$
$$=\frac{n}{n+1}\left(1+\frac{c}{X_{n}}-\frac{cn\E\left(\left.X_{n+1}-X_n\right|\F_n\right)}{X_n^2} +o\left(\frac{1}{n^{3/2}}\right)\right)$$
$$
\leq\frac{n}{n+1}\left(1+\frac{c}{X_{n}}-\frac{cg_{n}(X_n,0)}{(2+\beta)X_n} +o\left(\frac{1}{n^{3/2}}\right)\right)
$$
Recall that $g_n(x,0)=d-\frac{d(d-1)}{2}\frac{x}{(2+\beta)n}+o(\frac{x}{n})$ when $\frac{x}{n}\to 0$. Hence if $2+\beta=d$ we get
$$\E\left(\left.\frac{Y_{n+1}}{Y_n}\right|\F_n\right)
\leq 1-\frac{1}{n}+\frac{c}{X_n}-\frac{c}{X_n}\left(1-\frac{(d-1)X_n}{2dn}+o\left(\frac{X_n}{n}\right)\right)+o(1/n) 
$$
$$
=1-\frac{1}{n}+\frac{c(d-1)}{2dn}+o(1/n). 
$$
Hence, if $c<\frac{2d}{d-1}$ then $Y_n$ is supermartingale,  and by Doob's theorem there is a random variable $R$, such that $\sup_n Y_n<R$ almost surely. Hence for all $n$ 
$$X_n>\frac{cn}{\ln n + \ln R}$$
almost surely. Therefore 
$$\liminf \frac{\ln n X_n}{n}\geq \frac{2d}{d-1}$$
almost surely.
 


\section{Acknowledgment}

The present work was funded by a grant from the Russian Science Foundation (project No. 19-71-00043).


\begin{thebibliography}{Yury}


\bibitem[BA99]{BA99}
A. Barab{\'a}si, R. Albert.
\newblock Emergence of scaling in random networks.
\newblock {\em science}, 286(5439), 509--512, 1999.

\bibitem[BCDR07]{BCDR07}
C. Borgs, J. Chayes, Daskalakis, C. and   and Roch, S.
\newblock First to Market is not Everything: an Analysis of Preferential Attachment with Fitness.
\newblock STOC '07, June 11-13, San Diego, California, USA, 2007.

\bibitem[Chen03]{Chen03}
H.F. Chen. 
\newblock Stochastic Approximation and its Applications.
\newblock Nonconvex Optimization and its Applications, Springer, 64, 2002. -- 360 p. 




\bibitem[Gal16]{Gal16}
P.~A. {Galashin}.
\newblock {Existence of a persistent hub in the convex preferential attachment
  model}.
\newblock {\em PMS}, 36(1), 59--74, 2016.


\bibitem[HJ16]{HJ16}
J. Haslegrave, J. Jordan.
\newblock {Preferential attachment with choice.}
\newblock {\em  Random Structures and Algorithms}, 48, 751--766, 2016.

\bibitem[GLY19]{GLY19}
A. Grauer, L. L\"{u}chtrath, M. Yarrow.
\newblock Preferential attachment with location-based choice: Degree distribution in the noncondensation phase.
\newblock {\em ArXiv}, May 2019.
\newblock \url{https://arxiv.org/abs/1905.08481}.

\bibitem[HJY20]{HJY20}
J. Haslegrave, J. Jordan, M. Yarrow.
\newblock Condensation in preferential attachment models with location-based choice.
Random Structures and Algorithms, 56(3), 775--795, 2020.


\bibitem[Hof16]{Hof16}
R. Hofstag. 
\newblock Random Graphs and Complex Networks.
\newblock Cambridge University Press, Cambridge, 2016. -- 375 p.


\bibitem[JK77]{JK77}
N.~L. Johnson and S. Kotz.
\newblock Urn models and their application.
\newblock {\em John Wiley and Sons}, New York, 1977.

\bibitem[KR14]{KR14}
P.~L. {Krapivsky}, S.~{Redner}.
\newblock {Choice-Driven Phase Transition in Complex Networks}.
\newblock {\em Journal of Statistical Mechanics: Theory and Experiment}, P04021, 2014.


\bibitem[MP14]{MP14}
Y. Malyshkin, E. Paquette. 
{\it The power of choice combined with preferential attachement.}
Electron. Commun. Probab., 19(44), 1--13, 2014.

\bibitem[MP15]{MP15}
Y.~{Malyshkin}, E.~{Paquette}.
{\it The power of choice over preferential attachment}.
{ ALEA, Lat. Am. J. Probab. Math. Stat.}, 12(2), 903--915, 2015.

\bibitem[Mal18]{Mal18}
Y. Malyshkin.
{\it Preferential attachment combined with the random number of choices.}
Internet Math., 1--25, 2018.


\bibitem[Mah09]{M09}
H.~M. Mahmoud.
\newblock Polya urn models.
\newblock {\em Chapman and Hall/CRC}, 2009.

\bibitem[M{\'o}r02]{Mori02}
T.~F. M{\'o}ri.
\newblock On random trees.
\newblock {\em Studia Sci. Math. Hungar.}, 39, 143-155, 2002.

\bibitem[M{\'o}r05]{Mori05}
T.~F. M{\'o}ri.
\newblock The maximum degree of the {B}arab\'asi-{A}lbert random tree.
\newblock {\em Combin. Probab. Comput.}, 14(3), 339--348, 2005.


\bibitem[Pem07]{Pem07}
R. Pemantle.
\newblock A survey of random processes with reinforcement.
\newblock {\em Probab. Surv.}, 4, 1--79, 2007.



\end{thebibliography}
\end{document}